\documentclass[vecphys]{svmult}

\usepackage{makeidx}         
\usepackage{graphicx}        
\usepackage{multicol}        
\usepackage{graphics}
\usepackage{rotating}
\usepackage[bottom]{footmisc}

\makeindex             
\usepackage{times}
\usepackage{epsfig}
\usepackage{amsmath}
\usepackage{amssymb}

\newcommand{\bgo}{\begin{eqnarray*}}
\newcommand{\ego}{\end{eqnarray*}}
\newcommand{\bg}{\begin{eqnarray}}
\newcommand{\eg}{\end{eqnarray}}

\def\d3k{{\displaystyle {\rm d}{\bf k} \over \displaystyle (2\pi)^3}}

\begin{document}

\title*{Life and Times of Georgy Vorono\"{\i}\\ (1868-1908)}
\titlerunning{Georgy Vorono\"{\i}}
\author{Halyna Syta$^1$ \& Rien van de Weygaert$^2$}
\institute{$^1$ Drahomanov National Pedagogical University, 9 Pirohova Str, 01030 Kyiv, Ukraine\\
$^2$ Kapteyn Astronomical Institute, University of Groningen, the Netherlands\\
\texttt{syta@imath.kiev.ua; weygaert@astro.rug.nl}}

\maketitle


\printindex

\begin{verse}
\vskip -0.00truecm
\hfill{I know minutes not of complacency, not}\\
\hfill{of pride -- \ they all come later,  -- \ but moments}\\
\hfill{when the mind completely grips the idea which}\\
\hfill{before kept slipping off like a small ball.}\\
\hfill{Then I would forget about my existence}
\vskip -1.00truecm
\end{verse}
\begin{abstract}
Georgy Theodosiyovych Vorono\"{\i} (1868-1908) is famous for his seminal contributions to number theory, 
perhaps mostly those involving quadratic forms and Voronoi tessellations. He was born and grew up 
in the town of Zhuravka in the Ukraine, at the time part of the Russian Empire. Having studied at St.~Petersburg 
University under the supervision of Andrey Markov, in 1894 he became a professor of pure mathematics at 
the University of Warsaw. In his career he published six large memoirs and six short papers, each of which were so 
profound and significant that they left a deep trace in modern number theory. Together with Minkowski, he 
can be considered as the founder of the Geometry of Numbers. In this contribution, a brief sketch will be 
given of his life, work and legacy.
\vskip -0.5truecm
\end{abstract}
\section{Ukrainian Origins}
\begin{figure*}
\begin{center}
\vskip 1.0truecm
\mbox{\hskip -0.0truecm\includegraphics[width=11.5cm]{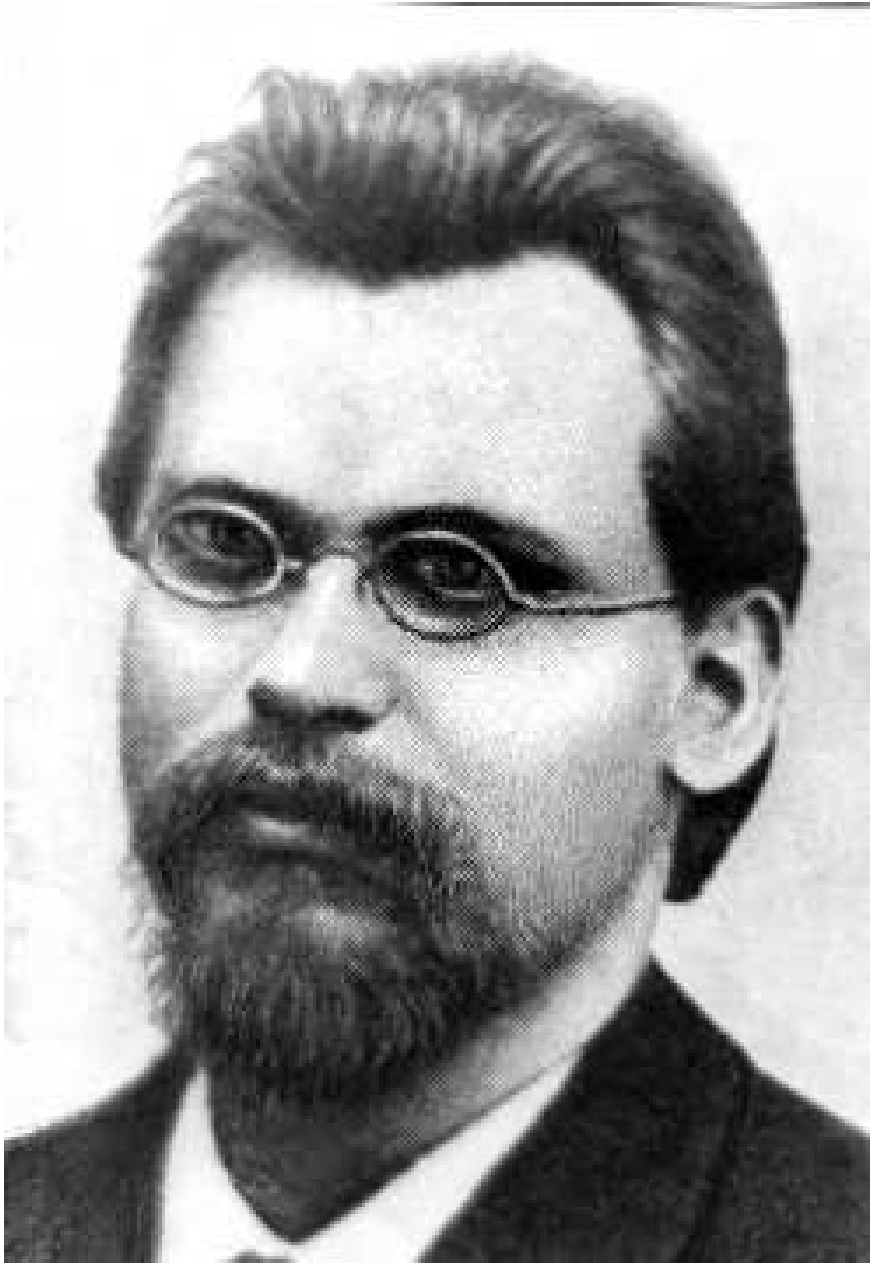}}
\end{center}
\begin{center}
\vskip 0.0truecm
\caption{Georgy Vorono\"{\i}, portrait.}
\end{center}
\vskip -1.0truecm
\label{fig:voronoi1}
\end{figure*}
\noindent Georgy Theodosiyovych Vorono\"{\i} \footnote{G.Vorono\"{\i} himself
used this transliteration of his name, that is "\"{\i}" at the end
of the word, in the papers written by him in French, whereas the
experts of Voronoi diagrams used to write "Voronoi" in accordance
with the term "Voronoi diagram". As it was accepted in Russia,
the second name "Theodosiyovych" is his patronymic. The first letter of 
his father's name was ``Fita'', the letter used in Russia in the nineteenth century, 
we reproduce it here as ``Th''. In Ukrainian
the name of the scientist is pronounced as: Heorhii Voronyi.} was
born on April 28, 1868,\footnote{According to the Julian calendar used
in Russia in the nineteenth century, in Russia this date was April, 16. Also the dates 
specifically mentioned by G. Vorono\"{\i} himself concern the Julian calendar. This 
concerns the dates mentioned on pp. 8, 9, 10, 11, 20, 21, 22 and 25.} in the small town 
of Zhuravka in the Poltava Gubernia, in Russia (now
the village Zhuravka belongs to the Varva District of the
Chernihiv Region in Ukraine). A picturesque hill over the Udai
river in Zhuravka appealled to one of Georgy's ancestors, perhaps his great-grandfather, 
who started out as {\it chumak} (an ox-cart driver; in medieval
times these carts were used for carrying the salt from
Crimea to Ukraine). After earning enough money, he bought a 
patch of land and settled there with his family \cite{Sh.P}.

According to the family legend, the Voronyis have acquired their
family name from a remote ancestor, a cossack captain  --
a Ukrainian military called {\it esaul} at the time of
Hetmanate, a commander of the Voronivka Fortress. There are
several settlements that carry such a name. On the map of De
Beauplan of the middle of seventeenth century, the fortified town of
Voronivka was on the right Dnipro bank, near Chyhyryn. And,
according to the description by M.~Tkachenko of the
Uman land on the right bank of the Dnipro river \cite{Tkach}, the town of
Vorone was near present-day Buky. There were mentions on Vorone
and of its fortress from 1545 to 1674, when the advance of
the Turks devastated much of the Uman area and forced the local
population to migrate en masse to the left bank. Perhaps some of
the mathematician's ancestors began as chumaks from the time of
this migration ...

\begin{figure*}[h]
\begin{center}
\vskip -0.25cm
\mbox{\hskip 0.0truecm\includegraphics[width=12.0cm]{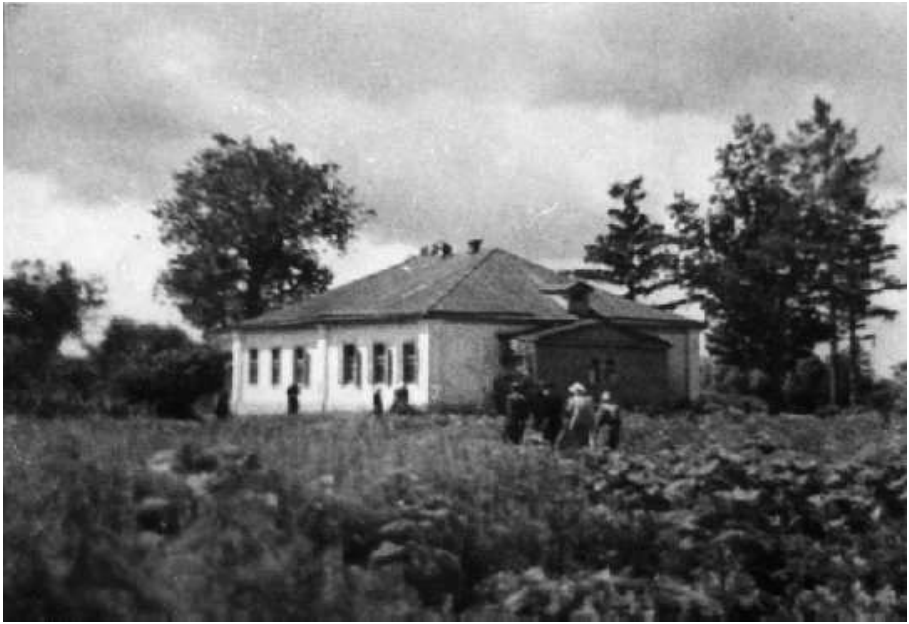}}
\end{center}
\begin{center}
\vskip -0.25truecm
\caption{The Vorono\"{\i} family house in Zhuravka. Picture of around 1953. 
Regretfully, the house does not exist anymore. The street at which the 
house stood is now called "Vorono\"{\i} Street".}
\end{center}
\vskip -0.8truecm
\label{fig:zhuravka}
\end{figure*}

Georgy's father, Theodosii Yakovych Voronyi (1837-1910), a son of ober-officer 
Yakov Tarasovych Voronyi, was educated at Kyiv University at the history 
and philology department (1857-1861). After his studies he worked at the Nemyriv gymnasium and 
the Nizhyn Lycee as a professor of Russian literature (1864-1872). Subsequently, 
he was a director of the gymnasium in Kyshyniv, then the gymnasium in Berdyansk and 
finally the gymnasium in Pryluky, some 18 kilometers from Zhuravka. Georgy was born 
while his father worked at  Nizhyn, and spent his years as a child in 
the towns of Nizhyn, Kyshyniv, Berdyansk, Pryluky, and Zhuravka. 
Theodosii was a man of progressive convictions and favoured efforts at 
popular education. Even in his student years at Kyiv University, Theodosii Voronyi
was involved in organizing free Sunday schools for working youth and
taught history in the Kyiv-Podil Sunday school
(1859-1861) \cite{Marakh}. He also hosted student literary
soirees. His initiative to create a Sunday school was supported by
the progressive circles amongst the Ukrainian intelligentsia. The famous poet
Taras Shevchenko \footnote{Taras Shevchenko (1814--1861) is a
great Ukrainian poet and artist. Many Ukrainian people worship him 
almost as their prophet. "Kobzar" is the famous book of poetry written by 
him.} visited this school in 1859 and  donated to the school 
50 copies of his "Kobzar". Later on, describing events of national
importance at Kyiv University in the late 1850s and early 1860s,
Olena Pchilka \footnote{ Olena Pchilka is a pen-name of Olena
Kosach Drahomanova (1849-1930), a well-known writer, ethnographer and a public
person. Her daughter, Larysa Kosach-Kvitka (1871-1913) became the 
famous Ukrainian poetess Lesia Ukrainka.} put a special emphasis
on the honorable actions of student Theodosii Voronyi \cite{Pchil}.

Theodosii Voronyi left us a manuscript of 15 pages (1861) in which he expressed his views 
on teaching and education at the Sunday School\footnote{It is kept at Kyiv City State archives, 
fund 16, discr. 471, 1861, file 105.}. In particular, he believed that "any
success in political and social life is impossible without people
being enlightened by moral sciences." He emphasized the need to
spread historical knowledge among the people. He believed that it 
"clarifies man's intellect, provides him with a better
understanding of his social status in life and frequently points
out the best ways of using one's abilities and to achieve
prosperity for oneself and wellbeing to others" \cite{Marakh}.

In 1887, Theodosii Voronyi retired. He got himself engaged in gardening in
Zhuravka and took an active part in the work of the Pryluky
Agricultural Society. Theodosii implemented his ideas of popular
education in Zhuravka. At his own expense he built a school
for the village children. The school also became the place where
lectures, concerts and performances were held. The money collected
was spent to enlarge the public library in Zhuravka, opened with
the help of the Voronyi family.

Theodosii had a daughter and three sons. His elder son, Leonid,
became a doctor, his son Mykhailo was an agronomist. Georgy's fate was to
become a well-known scholar. 

Mykhailo inherited his father's passion for gardening and took
part in creating a small company for cultivation and herb
processing in Zhuravka (1885). This was one of the first
companies (if not even the first) of such kind in Russia. Thanks
to the Voronyi family, Zhuravka played a role as a kind of regional 
cultural centre.

\begin{figure*}
\begin{center}
\mbox{\hskip 0.0truecm\includegraphics[width=8.0cm]{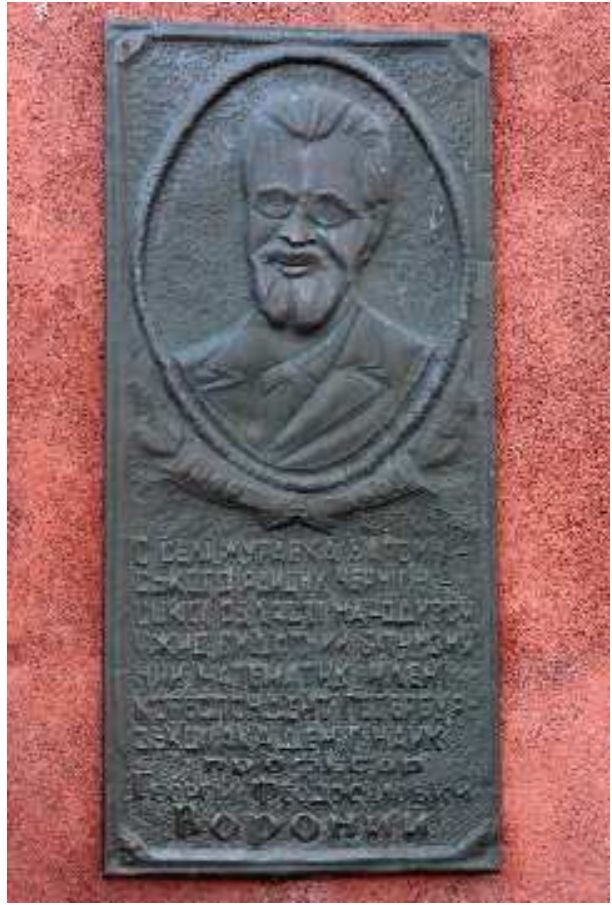}}
\end{center}
\begin{center}
\vskip 0.0truecm
\caption{The Vorono\"{\i} school in Zhuravka: the commemoration plaquette of Vorono{\"\i} at the 
entrance wall was made by local inhabitants. The school was built next to the location of 
Vorono\"{\i}'s old elementary school. It has a room with exhibition on 
Georgy Vorono{\"{\i}}.}
\vskip -1.0truecm
\end{center}
\label{fig:school}
\end{figure*}

\section{Schooldays: Zhuravka \& Pryluky}
\noindent Georgy Vorono\"{\i} first followed school at Berdyansk, and later at the
Pryluky gymnasium from which he graduated in 1885. While a student at the Pryluky gymnasium, 
Georgy already distinguished himself by his deep interest in science, and by his diligence 
and punctuality. 

His pupil's testimonial report stated that ".. with excellent talents, in 
spite of his young age, he has reached a considerable mental maturity and 
displays a very serious interest and love for learning. He acquired himself 
a very good knowledge in all subjects taught at high school. He has a 
particular inclination and calling for mathematics, in which he has 
considerably surpassed the average student's progress'' \cite{Sh.P}. 

He was particularly attracted to mathematics, which he studied 
very thoroughly and in which he revealed his brilliance, being 
especially fond of algebra. Georgy's interest in mathematics was stimulated 
by his favourite teacher of mathematics, Ivan Volodymyrovych Bogoslovskii, who 
greatly influenced his outlook and attitude towards the world. Later on, in his 
student years, Vorono\"{\i} appealed several times to his beloved teacher and his 
father as his moral and mental guides. Characteristically, in these years Georgy 
also kept evaluating himself again and again to see if he would have enough talent 
and persistence to be a professional scientist. 

\begin{figure*}[h]
\begin{center}
\vskip -0.25truecm
\mbox{\hskip -0.0truecm\includegraphics[width=11.5cm]{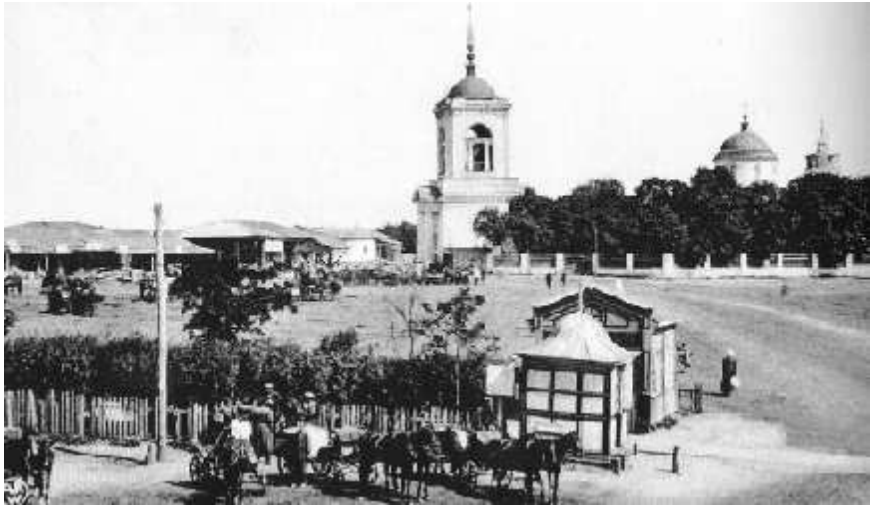}}
\end{center}
\begin{center}
\vskip -0.25truecm
\caption{Pryluky, view of the market place. The picture has been taken around the time 
Georgy Vorono\"{\i} attended the gymnasium of Pryluky, around the year 1885.}
\end{center}
\vskip -1.0truecm
\label{fig:priluky}
\end{figure*}

While at gymnasium, he solved a problem posed by Professor Ermakov of 
Kiev University on factorising polynomials, which had appeared in the Journal of 
Elementary Mathematics (Kyiv). In 1885 the young Georgy Vorono\"{\i} presented his 
solution in the paper {\it "Decomposition of the polynomials on factors based 
on the properties of the roots of quadratic equation"} \cite{vor1885}. Although still 
a school pupil, this led to his first publication, in the same journal. It got published 
in the same year he entered St.~Petersburg University.

\section{University Student: St.~Petersburg}
\noindent In 1885 Georgy Vorono\"{\i} entered the Faculty of Physics and Mathematics 
of the University of St.~Petersburg, where he studied mathematics during the period 
1885-1889. For mathematical education this was at that time the best university in 
Russia: instruction in mathematical disciplines was being given there at the high level 
to which it had been raised by the genius of P.L.~Chebyshev and the staff of the celebrated 
mathematical school of which he was the head (also see \cite{sako1968}).

At the university Vorono\"{\i} attended the lectures of A.A.~Markov, A.N.~Korkin and 
Yu.V.~Sokhotsky.  He was particularly attracted to algebra and found Sokhotsky's course 
the most enjoyable of all. In his diary he wrote: ``The pure mathematics lectures 
captivate me more and more. I prefer Professor Sokhotsky's lectures in the special 
course on higher algebra to all the others. ... The main thing that concerns me is 
whether I have enough talent.''

There is no doubt that Vorono\"{\i} had ample talent and his subsequent 
accomplishments show this very clearly. Nonetheless, his student years were very difficult for 
him. In particular financially Georgy did not have an easy time. Even though his father had held 
posts of considerable standing, while retired he was not able to fully support his son. Georgy 
had to earn his living by private tutoring of mathematics, often for rather 
small fees. Taking such work very seriously and putting much effort 
into it, he found giving lessons exhausting \cite{conno07}. As he wished to spend 
as much time as possible on his own investigations and mathematical research, he had 
to restrict the tutoring to the level of providing him with the bare essentials. 
He also found life in a student residence difficult, as he had to live in a cold, 
damp room. In his diary he frequently writes about the difficulties in accustoming 
himself to working under straining conditions. The constant amount of loud noise 
made by his roommates, the frequent student disturbances, and working at 
night did not make his life easy. It was only when he was in his last year, that 
Vorono\"{\i} was offered a very meagre student bursary. 

\begin{figure*}[h]
\begin{center}
\mbox{\hskip 0.0truecm\includegraphics[width=11.5cm]{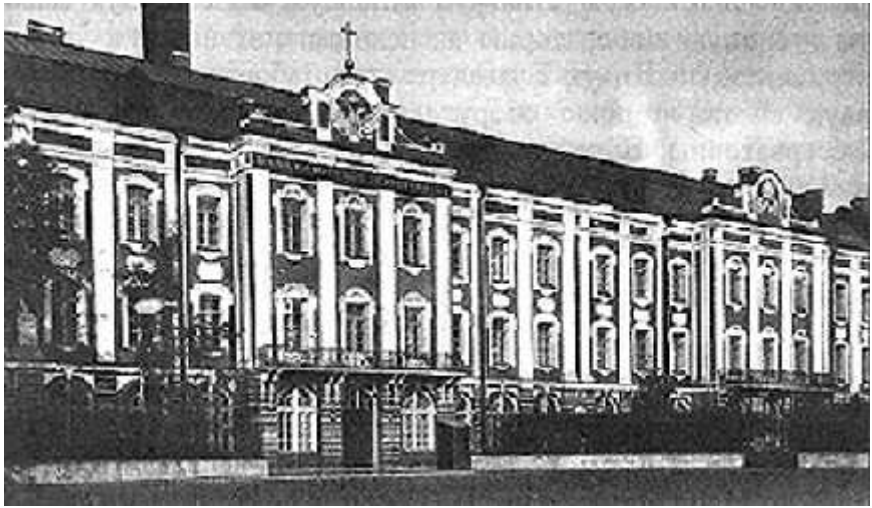}}
\end{center}
\begin{center}
\vskip 0.0truecm
\caption{The central building of St.~Petersburg University, the Twelve Collegium, at the 
end of the 19th century.}
\end{center}
\vskip -1.0truecm
\label{fig:university}
\end{figure*}
\noindent All these sad self-revelations do not prevent Georgy from studying his 
favorite science with fervour and perseverance. He developed a
strong self-discipline, marked by the habit of regular and diligent 
work and the ability of concentrating his intellect on the problems 
at hand.  

\subsection{Georgy's student diary}
It says something about the personality of Georgy Vorono\"{\i} that in 
these student years he confided his doubts to his diary. Fortunately, this diary 
has been partially preserved. Along with his descriptions of everyday experiences 
and events, it is a sincere self-confession of a young man. It
discloses his character, his inner world, the process of his 
creative growth and self-consciousness. The author is active and
sensitive and cannot remain indifferent to the events around him. 
He also tries to help when necessary. At times he is hot-tempered, for 
which he later expresses regret. He states "I am merrily gazing at 
God's world and to everything I touch I submit myself with rapture". Georgy 
aims "to reach everything by heart, and not just by intellect" 
and tries to look at himself from the outside. In this, he displays a 
rather low self-esteem, while also trying to grasp his own feelings and 
inclinations:

\begin{figure*}[h]
\vskip -0.25truecm
\begin{center}
\framebox[8.5cm]{\hskip 0.0truecm\includegraphics[width=8.0cm]{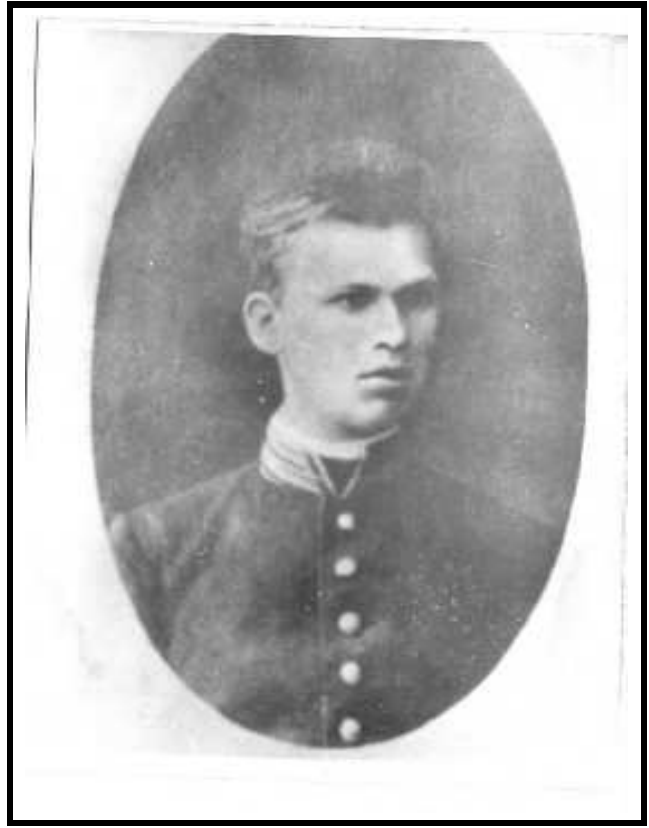}}
\end{center}
\begin{center}
\vskip -0.25truecm
\caption{Georgy Vorono\"{\i} in his student years, 1885-1889.}
\end{center}
\vskip -1.0truecm
\label{fig:student}
\end{figure*}

"What am I after all? I am fond of playing cards. I do not have 
any noble pride. That is, if I am mocked I do not get angry and
do not quarrel with the offender. I feel my weakness in front of the
powerful of this world"...
\smallskip

The author has a rather melancholic view of contemporary
society: "Our times are hard, we are victims of a terrible regime. The most 
innocent things cannot be said. Otherwise one gets into the hands of the 
{\it custodians of hearts}. We are characterized by mistrust. We distrust each 
other." The false feelings and hypocrisy that reigned around him led him to
seclude himself from social life. 

We may obtain some insight into his lifestyle during these student days 
from remarks like those on his study of astronomy: "For the second day I have
been busy calculating the solar eclipse of August 7. Yesterday I was
working for 10 hours, today for 7 hours, the work has progressed 
considerably but I am feeling quite exhausted, the more so
because I have not had a breath of fresh air for two days.
Figures, figures ... Yesterday I so crammed my head with them,
that they pestered me all the night and I had to get up and 
throw cold water on my head. Maybe today I shall have to resort 
to the same means." Ten days later he wrote in his diary: 
"I have not been writing my diary for ten days. It was not because of 
being in a lazy mood, but just because I am short of time. All the Holy Week 
I was busy at calculating the solar and lunar eclipses. Concerning the former 
I seem to obtain a wrong result. So I'll have to redo some parts of the 
calculation."

In spite of his studiousness, Georgy failed the astronomy exam. 
It happened for the first time in his life. "This failure
made a terrible impression on me, but I hid it; as usual it
provoked in my mind a lot of questions: am I really a failure even
here~? No, I shall not agree with it. I can call to witness my
comrades, they will confirm that I possess enough knowledge. But
the main thing was that it all awoke in me a passionate desire
for study".

\section{Research Student}
The main subject that Georgy Vorono\"{\i} chose for his specialization was 
number theory. In the second half of the 19th century, the leading
professors of St.~Petersburg University, -- P.L.~Chebyshev, Yu.V.~Sokhotsky, 
A.A.~Markov, Ye.I.~Zolotaryov, -- actively worked in this field. The mathematical school 
later called the St.~Petersburg School of Number Theory arose as 
a result of their activities. It was Professor Andrey Markov who 
became Georgy Vorono\"{\i}'s main mentor and supervisor \footnote{Andrey 
Andreyevich Markov (1856-1922) was a famous Russian mathematician. He is 
best known for his work on the theory of stochastic processes, 
resulting in well-known and important concepts such as Markov 
chains and Markov process.}.

\subsection{Scientific debut: Bernoulli numbers and the Staudt Theorem}
Georgy made his debut in the mathematical circle on December, 2, 1888 
with a presentation on his own investigation on some properties of 
Bernoulli numbers. Professor Markov received Vorono\"{\i}'s
presentation with great approval. This first success increased
Georgy Vorono\"{\i}'s enthusiasm, encouraging him to continue and persist 
in his research work with even greater intensity. 

In a timespan of only three weeks Vorono\"{\i} succeeded in proving 
the assertion stated earlier by Adams. The latter had done so without proof, 
posing the following proposition in his paper, "If $p$ is a prime
divisor of $n$ and it is not a factor of the denominator of the
$n$-th Bernoulli number, then the numerator of $B_n$ is divisible by
$p$" \cite{adams1878}. In this connection Adams had 
noted: "I have not succeeded however, in obtaining a general proof of 
this proposition, though I have no doubt of its truth".

Vorono\"{\i} continued his research and by following the same methods, 
he obtained the next result a couple of days later: the Staudt theorem. He 
also found several new generalizations. He handed his new results to professor 
Markov. 

\begin{figure*}[h]
\begin{center}
\mbox{\hskip 0.0truecm\includegraphics[width=8.0cm]{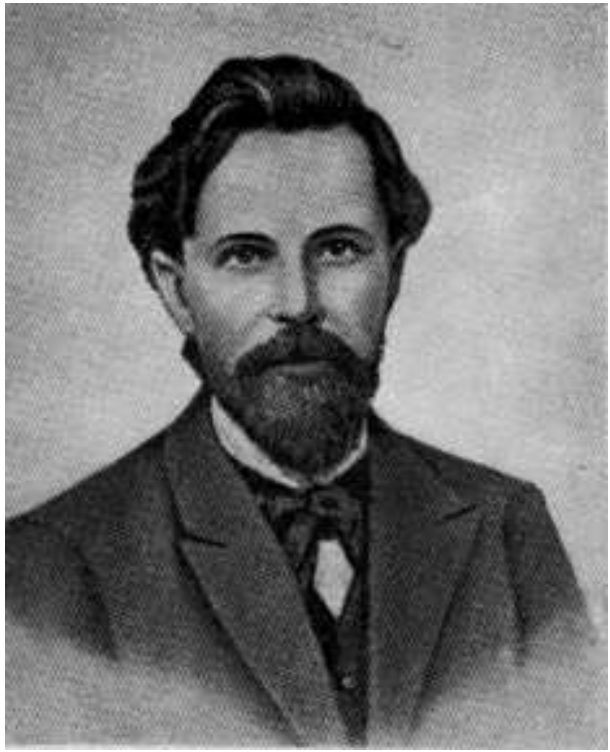}}
\end{center}
\begin{center}
\vskip 0.0truecm
\caption{Andrey Markov (1856-1922), master and doctoral thesis supervisor of Georgy Vorono\"{\i}. Markov 
is famous for his contributions to the theory of stochastic processes, resulting in well-known 
and important concepts such as Markov chains and Markov processes.}
\end{center}
\vskip -1.0truecm
\label{fig:markov}
\end{figure*}

On the eve of the new year 1889, Georgy wrote a brief report and evaluation of 
the past year, as he had become accustomed to do: 

\medskip
"Well, this year has not passed in vain for me. I have been working a lot, 
very much, and made certain that I can work and, it seems to me, succeeded 
in making others to be convinced of this.

At this time last year I had expressed a timid wish in the diary and 
now I see  it has come true. The thing, which I was afraid of, does not exist. 
I know, I do believe, that on the basis of science work, and only that,  
I shall find my good luck [...] I am not a poet and  I do not know the 
inspiration which poets feel, but I know minutes not of complacency, not of 
pride -- they all come later -- but moments when the mind completely grips
the idea which before kept slipping off like a small ball. Then I
forget about my existence. 

I firmly believe that in this respect the new year will bring me still more joy, 
because I noticed that for my latest successes I am obliged to a habit of thinking 
without a pen and paper. All assertions, which I proved, occurred to me quite
independently, and I only had to verify them. I hope this habit to
think in such way will stand me in good stead".
\medskip

On January 9, 1889 Georgy got a positive response from Markov on
his new results. Markov proposed to him to prepare the paper for
publication, but Georgy Vorono\"{\i}, demanding to himself, first 
wanted to refine the result obtained and to generalize it. Markov
told that he had looked through the table of Bernoulli numbers and
doubted whether it would be possible to find the law of the 
numerator, because it involves a very large number of simple quantities. 
In response (as Georgy wrote in the diary) Vorono\"{\i} presented to him 
the theorem which introduced these quantities. 

\medskip
"Markov was greatly interested by it, as the theorem which I
had presented to him before gave nothing of the kind. I told him I had 
so far no full proof of that theorem, but that I would 
hope to prove it soon...

I came home feeling quite jaded, unwilling to do anything, and I
felt at this time the whole burden of solitude. But the habit
insisted on its own, I went to bed soon and today from 6 o'clock
I have been working to prove that proposition. If anybody asked
me where I had taken that theorem from I should certainly find
difficulty in replying. I did not invent it, only proved it, but in
leaps, and filled gaps between them by imagination. But then I
checked it on very big examples and never got a contradicting
result. Even today when I felt some doubt I checked its
correctness with the help of the 44th Bernoulli number, which
has 70 figures in the numerator. Though probability of a mistake
was only $1/17$, just on the basis of selected conditions I got 
to doubt. Again I got a confirmation of my theorem. Had there been
any misprint in the 70-figure number, it would have grieved me so
deeply, but fortunately it did not happen."
\medskip

Georgy was so much carried away by his new research that he 
could not make himself do anything else. For example, to
prepare for the exams...

\medskip
"An idea flashed across my mind: why shouldn't I occupy myself
with some problem in mathematics -- it seemed as if a ray of light fell
on me, but then went out at once: I cannot, no time (even though I do nothing
at all)... Yesterday I wrote down these lines, but once I sat to work 
within a few hours I revised all my former work. 
All that was vague may be understood in such a simple way
that one cannot have wished for a better one."
\medskip

That year Vorono\"{\i} continued to work persistently during his
summer vacation in Zhuravka. Only when his exposition became
sufficiently clear and precise he sent the paper, {\it ``On Bernoulli numbers''} 
to the Communications of the Kharkov Mathematical Society \cite{vor1890}. In the autumn 
of 1889, Georgy Vorono\"{\i} passed his final exams brilliantly and defended his 
candidate thesis\footnote{The term "candidate thesis" corresponds to the present term 
"degree work".}. The theme was Bernoulli numbers. 

In November 1889 he stayed at the University to prepare for his candidate degree exams.
The representation necessary for his degree was signed by all the leading
mathematics professors: A.A.~Markov, O.N.~Korkin, Yu.V.~Sokhotsky and 
K.A.~Posse. At the suggestion of these he was allowed to remain at the university 
for further studies leading to the title of professor. While gaining a scholarship to 
enable the continuation of his studies, Georgy Vorono\"{\i} was also appointed as a
supernumerary teacher at the Peterhof Progymnasium. From 1889 to 1894 he worked at this school 
as a civilian instructor in mathematics. 

\begin{figure*}[h]
\begin{center}
\vskip 0.0truecm
\framebox[11.8truecm]{\hskip 0.0truecm\includegraphics[width=11.5cm]{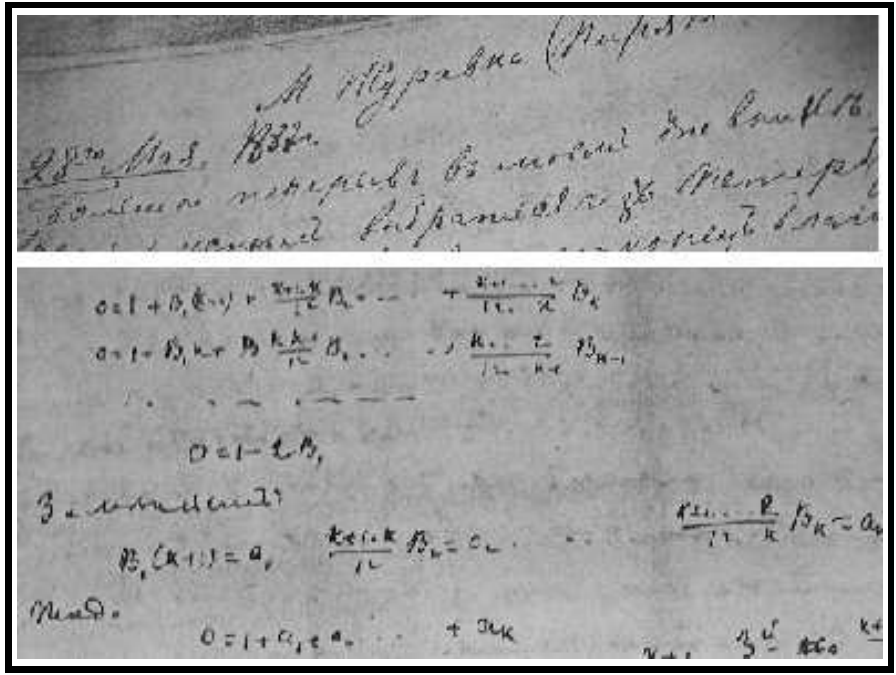}}
\end{center}
\begin{center}
\vskip -0.2truecm
\caption{Excerpts from the diary of Georgy Vorono\"{\i}. Just below the centre of the image we 
can distinguish the phrase {\it Zhuravka of Pyryatyn uyezd of Polatava gubernia} (``uyezd'' means 
district, in tsarist Russia ``gubernia'' was the word for province).}.
\end{center}
\vskip -1.0truecm
\label{fig:diary}
\end{figure*}

\subsection{New Year 1890: promise and prospects}
\noindent On the eve of the new year 1890 Georgy, following his habit, sums
up his accomplishments of the past year:

\noindent ".. Everything that I did not even dare to dream of came true: I
passed the exams, having been awarded fives for all of them. I could stay at 
the University with the help of the grant. In a word, to a great extent my future has 
already taken shape, even though it happened at the last moment. Until then I had only been 
suffering, enduring my pupil, who vexed me with his laziness. When I came home, I fell 
exhausted on the bed and thought: when will there be an end to it? Alas! I had to get 
up and swot, swot on and on! 

If I add that at this time I had influenza, was without a copeck
of money, had no prospect except for the awful private lessons,
was living in a damp cold room, which I remember now as a
nightmare, it is strange indeed that I did not collapse.  

Passion for research, for finding new properties and relations between 
quantities has developed in me to an inconceivable extent; I can hardly put down my pen. 
The most urgent things and obligations -- everything steps to the background, 
and I go on writing and writing...

I often compare myself with an alchemist, because, like him, I
have no guiding star and have only a passion. And this passion has
developed to such extent that I am losing my sleep as soon as it
seems to me I have touched anything of importance. But alas, so
far after a sleepless night I could see only that I had run
against a solid wall and just was nourishing my illusions.

I am not embarassed by it at all, as I have already become firmly
confident that I can easily take a simpler problem and get some
result -- not an essential one, only formal; but it does not tempt
me, these laurels I have already reaped and now can say like
Themistocles: "Gauss' laurels keep me awake!"

\section{Olia}
\noindent Georgy usually spent his vacations in Zhuravka, in his 
"native palestines", as he wrote. He visited the nearby village of Bohdany
quite often to see his beloved girl Olia Krytska, his future wife.

Recollections about his acquaintance and the development of his
relations with Olia Krytska occupy a particular place in the
diary. Georgy writes so sincerely about his feelings, with such
virtue and temperament -- (\/events are almost ignored, only his
feelings are recorded\/) -- that these pages read like a real
novel.  He determined once and forever for himself that his
destiny was in Bohdany, but he concealed his feelings for the time
being because he had no financial basis for his own family. His
father insisted on this decision. Such a vagueness in relations
brought him many sufferings, but he patiently waited for his hour
and did not permit any other passion to find the way to his heart.
In 1889, on the eve of his departure, Georgy wrote about his last visit to
Bohdany:

"Once more I am writing down my last visit to Krytskis... 
I am mounting the horse, once more saying goodbye to everybody, 
that is the end to everything which filled my life during the
four months and which will cause me to behave stern and cool
during the whole stretch of the Petersburg year.

Only mathematics as a bright star is shining afore me, in it I
trust all my hopes... The experience of the last year has strengthened 
my endurance, and my creative eagerness, suppressed before, is bursting 
into action, and I am certain that Petersburg will bring me much that is 
new in this respect.

So goodbye, Olia, goodbye, Zhuravka! Till the new spring I shall
cover myself with my armour. And, as if dreaming I shall see  this
summer, which gave me so much strength and health and those grains
of happiness, which I know I shall so often experience when
reading my diary in Petersburg, picking them from those talks
with Olia, which I wrote down, along with everything which so often 
made my heart beat."

\begin{figure*}[h]
\begin{center}
\vskip -0.25truecm
\mbox{\hskip -0.18truecm\includegraphics[width=12.4cm]{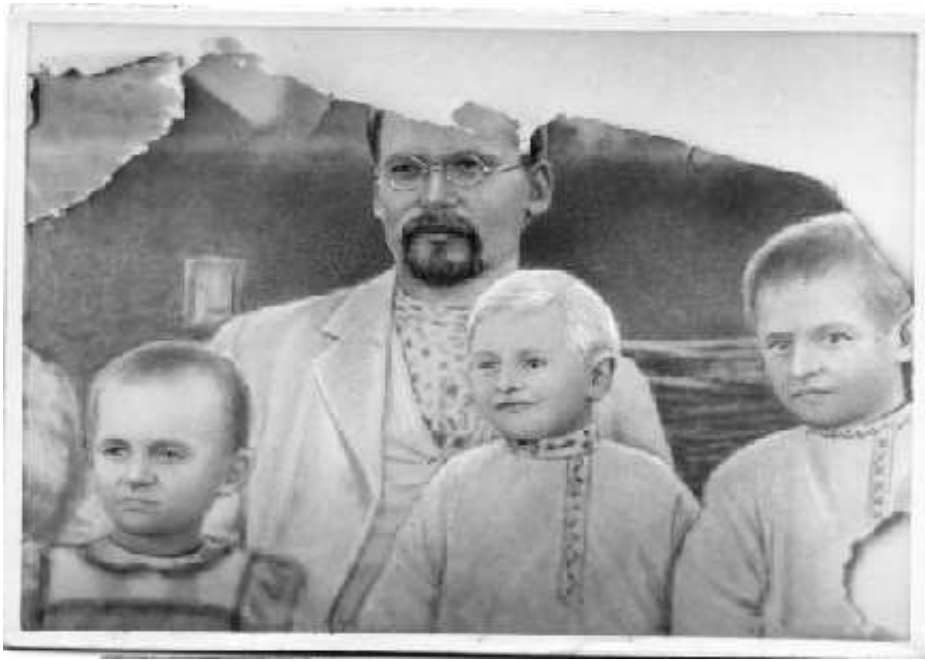}}
\end{center}
\begin{center}
\vskip -0.25truecm
\caption{Georgy Vorono\"{\i} with his three eldest children,
Oleksandra, Yuri and Oleksander. Zhuravka, 1904-1905.}
\vskip -1.0truecm
\end{center}
\label{fig:family1}
\end{figure*}

\begin{figure*}[t]
\begin{center}
\vskip -0.1truecm
\mbox{\hskip -0.18truecm\includegraphics[width=11.4cm]{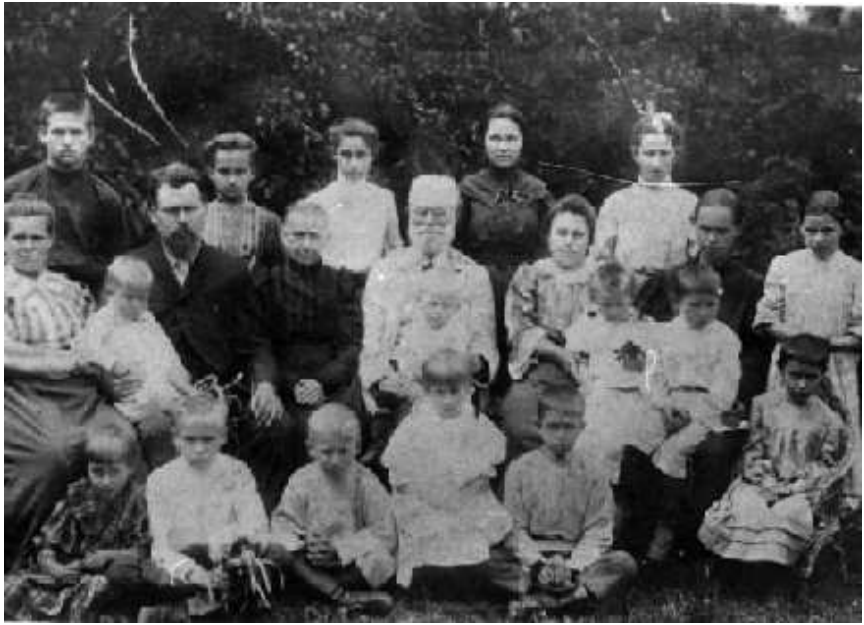}}
\end{center}
\begin{center}
\vskip 0.0truecm
\caption{Photo of the large Georgy Vorono\"{\i} family (approx. 1905-06, Zhuravka). 
Central row, sitting from left to right: Olha Mytrophanivna Vorona Krytska (Georgy 
Vorono\"{\i}'s wife); Georgy Vorono\"{\i} with his daughter Tetyana (she died
in  childhood); Kleopatra Mykhailivna Vorona Lychkova (Georgy Vorono\"{\i}'s mother);
Theodosii Yakovych  Voronyi (Georgy Vorono\"{\i}'s father), with his grandson 
Igor (son of Mykhailo Voronyi, Georgy Vorono\"{\i}'s brother); Evfrosiniya Ivanivna 
Vorona (Mykhailo Voronyi's wife) with her daughter Halyna; Nadiya Theodosiivna Ermakova 
Vorona (Georgy Vorono\"{\i}'s sister) with her son Peter -- Nadiya became a widow quite 
early, and had a large family of seven children, she was helped by Georgy Vorono\"{\i} 
in her family problems; Vira Petrivna Prosvirlina-Ermakova (N.Th.Ermakova's daughter). Other 
persons around are their children.}
\end{center}
\vskip -0.8truecm
\label{fig:family2}
\end{figure*}

As Georgy Vorono\"{\i} gained his scholarship following his candidate's degree in 1889, and 
got appointed as supernumerary teacher at the Peterhof Progymnasia, his future 
was finally ensured. At last, he became Olia's fianc\'e. On the eve of 1891, 
he wrote in his diary:\\

\bigskip
{\parbox[c]{10.5truecm}{\small
\hfill December 31, 1890, Peterhof\\
\ \\
\ \ \ True to the old custom, today, on the eve of New Year, I cast a glance at how 
I have lived through and deeply felt the Old Year. The first thing which I gladly  
note and which has become a harbinger of my future happiness is: {\it Olia loves 
me}. I know it now for certain ! How happy I am ! So long I had been silently 
suffering from doubts, and at last it has been clarified, and I have already 
become Olia's fianc\'e ! ... }}

{\parbox[c]{10.5truecm}{\small
\vskip 0.5truecm
\ \ \ Yes, now I know well that Olia loves me, but nevertheless lasting doubts and 
expectations have brought some bitterness. I seem to have become hardened in my 
permanent solitude. Ever growing passion for Mathematics has developed in me an 
egotism of no small degree. I am afraid I cannot feel strongly and surrender fully  
to my feelings. 

\ \ \ As for me the mind comes ahead always and everywhere. And the worldly wisdom, 
known from books, is saying that mind and love can scarcely be reconciled. That is 
what makes me fear sometimes that Olia probably will not be happy with me. As for me, 
I shall probably always take refuge in Mathematics.}}
\bigskip
\medskip

\noindent Georgy and Olha married in 1891. They got six children, two sons and four daughters: 
Oleksander (1892), Oleksandra (1894), Yuri (1895), Maria (1900), Tetyana (1904) and 
Olena (1906). One of his daughters, Tetyana, died in childhood. 

\begin{figure*}[h]
\begin{center}
\mbox{\hskip 0.0truecm\includegraphics[width=12.0cm]{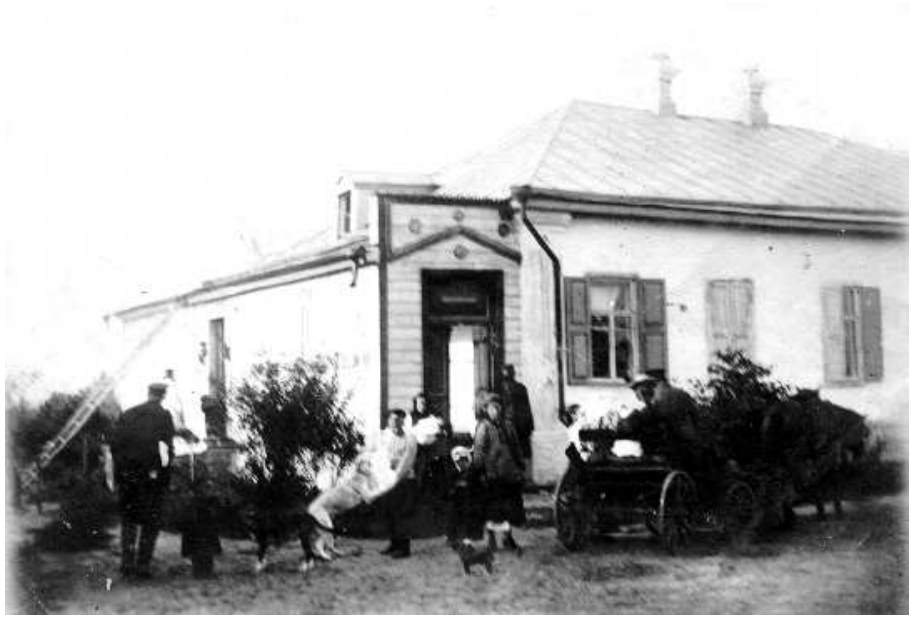}}
\end{center}
\begin{center}
\vskip 0.0truecm
\caption{End of summer: Vorono\"{\i} and his family leave the Zhuravka house for 
Warsaw, to take up his post at Warsaw University.}
\end{center}
\vskip -1.00truecm
\label{fig:zhuravkawarsaw}
\end{figure*}

His family followed him to his posts. While throughout his short life he kept  
spending his summers in Zhuravka, Olha and the children moved back with him to 
St.~Petersburg and Warsaw to stay with him throughout the rest of the year. 

\section{Master and Doctoral Thesis: the Vorono\"{\i} Algorithm}
Having decided to stay at St.~Petersburg University, the scientific interests of the young 
researcher concentrated more and more on the theory of irrationalities of the third degree. 
He wrote two large studies on this subject, one which formed the basis of his 
master dissertation and the second one being his doctoral thesis work. Both Vorono\"{\i}'s 
master thesis and doctoral thesis were of such high quality that they were awarded the 
Bunyakovsky prize by the St.~Petersburg Academy of Sciences \cite{conno07}.

\subsection{Master Thesis: integral algebraic numbers}
His master's dissertation concerned the issue of algebraic integers associated with the 
roots of an irreducible cubic equation, ``On algebraic integers dependent on a root  
of an equation of the third degree''. He defended the master thesis in 1894. In the 
preface he wrote \cite{conno07}:

``In the essay I am now presenting, results from the general theory of algebraic integers are 
applied to the particular case of numbers depending on the root of an irreducible equation $x^3=rx+s$. 
The results obtained turn out to be very graphic. All integers of the field in question have a certain 
form. Using the form taken by the integers, it is not difficult to find a form embracing all the integral 
numbers divisible by a given ideal number, or, in other words, to find the ideal corresponding to that 
given ideal number. In our exposition the resolution of these questions is based on a detailed study 
of the solutions of third-degree equations relative to a prime and a composite modulus.''

\subsection{Doctoral Thesis: the Vorono\"{\i} Algorithm} 
The other study concerned his doctoral thesis, titled ``On a generalization of the algorithm for continued 
fractions'', which was published in 1896 and which he defended in 1897. 

With the help of continued fractions it became possible to approximate the quadratic irrational numbers,  
because each of these numbers may be represented by the periodic continued fraction. Generalization of 
this algorithm to the cubic irrational numbers appeared to be so difficult that it could not be solved 
during the nineteenth century. G.F. Vorono\"{\i} was the first mathematician who managed to find a 
generalization. 

It took another 42 years until the Vorono\"{\i} algorithm was rediscovered by G. Bullig \cite{Bull}. In 1976 
I.O. Angell \cite{Ang} used it in his tabulation of totally real cubic fields with a positive discriminant 
$\le 100000$ \cite{Schi}. A generalization of the Vorono\"{\i} algorithm has been proposed by J. Buchmann 
in 1985 \cite{Buch1,Buch2}. 

\medskip
\noindent In the preface to his thesis Vorono\"{\i} writes:

``In his article 'De relatione inter ternas pluresve quantitates instituenda' Euler gave the first generalisation 
of the continued fraction algorithm. ... Jacobi changed Euler's algorithm somewhat in order to unify the 
calculations''. Vorono\"{\i} then discussed Jacobi's generalisation, after which he went on to look at generalisations 
due to Poincar\'e and Hurwitz. He then looked at contributions by Dirichlet and Hermite, showing that none of the 
above provided a satisfactory generalisation. He then writes: ``Taking as the basis of our investigations a special 
way of viewing the continued fraction algorithm, we propose a new generalisation of it.'' \cite{conno07}.

\medskip
Within the context of this work, Georgy Vorono{\"{\i}} addressed the issue of the analysis of cubic fields 
and the ring of integer algebraic numbers of such fields. The units of the ring are the invertible elements 
$a \in K$ for which there exist elements $b \in K$ such that $ab=1$. The result obtained was so striking that 
prof. Andrey Markov could not believe the correctness of Vorono\"{\i} proofs and did not
dare to approve his work. In this connection  D.~Grave asserted
\cite{Gr.D}:

\medskip
"Markov asked Vorono\"{\i} by telegraph to come from Warsaw to Petrograd
(this fact is known to me from Markov himself) \footnote{At the time of these events 
the city was called St.~Petersburg. Grave used here the name Petrograd, in accordance
with the city's name in the 1930s.}. Markov invited Vorono\"{\i} to his office and proposed him 
to calculate the unit for the cubic equation $r^3 = 23$. By artificial means, Markov had found for this 
example the unit
\begin{equation}
e = 2166673601+761875860r+267901370r^2\,.\nonumber
\end{equation}
\noindent Vorono\"{\i} calculated for three hours. The period had 21 terms and in
order to find the main unit it was necessary to multiply 21 expressions
\medskip
\begin{equation}
-2+\rho,\ {-{11}+2\rho+\rho^2\over{15}},\
{-3-\rho+\rho^2\over4},\ {-9+5\rho+\rho^2\over{17}},\
{4-3\rho+\rho^2\over{10}},\nonumber 
\end{equation}
\begin{equation}
{1-\rho+\rho^2\over8},\
{-2+\rho\over3},\ {1+3\rho-\rho^2\over{10}},\
{-5-\rho+\rho^2\over3},\ {-1+\rho\over2},\nonumber
\end{equation}
\begin{equation}
{-{10}+\rho+\rho^2\over{11}},\ {-2+\rho},\
{-{11}+2\rho+\rho^2\over{15}},\ {1-\rho+\rho^2\over8},\
{-2+\rho\over3},\nonumber
\end{equation}
\begin{equation}
{1+3\rho-\rho^3\over5},\ {-1+\rho\over2},\
{-1+{10}\rho-\rho^2\over{33}},\ {-{11}+7\rho+\rho^2\over{20}},\nonumber
\end{equation}
\begin{equation}
{9-7\rho+2\rho\over{31}}, \ {-5-\rho+\rho^2\over6}.\ \nonumber \\
\ \nonumber 
\end{equation}
\medskip
\noindent Following this analysis, he found the unit 
\begin{equation}
E = -41399 - 3160r + 6230r^2\,.\nonumber
\end{equation}
\noindent It turned out that $Ee = 1$. So, it was verified that the algorithm
really worked."  \cite{Gr}.

\medskip
In connection with his doctoral thesis, the academician D.~Grave wrote in 1934 \cite{Gr.D}:
\ "Georgy Vorono\"{\i} is a Ukrainian mathematician of genius.
While being at Petersburg University he was busy with an 
investigation of the cubic sphere which he performed with remarkable
success leading to a great discovery in this field. He extended
the algorithm of continued fractions for the cubic area which gives
algebraic unities in the quadratic area. This extention had been
sought in vain by all of the most outstanding mathematicians throughout 
the 19th century. This is how the Vorono\"{\i} algorithm was found".
\medskip

\section{Professor in Warsaw}
\noindent After defence of his master dissertation in 1894, Vorono\"{\i} was
appointed as professor of pure mathematics at the Russian Imperial University of Warsaw 
(today's Warsaw University), where he worked almost all the rest of his life. While continuing his scientific 
research on his doctoral thesis work, his lecturer's duties took too much time 
because the staff of the University included, in addition to Vorono\"{\i}, only two 
professors. For this reason, Vorono\"{\i} was forced to lecture
for the  students of different terms gathering them in one class.
For example, Vorono\"{\i} lectured the courses of Number
Theory  and the Probability Theory to students of the third and
fourth years. These courses were given by Vorono\"{\i} once every two
years.

\begin{figure*}[h]
\begin{center}
\mbox{\hskip 0.0truecm\includegraphics[width=8.0cm]{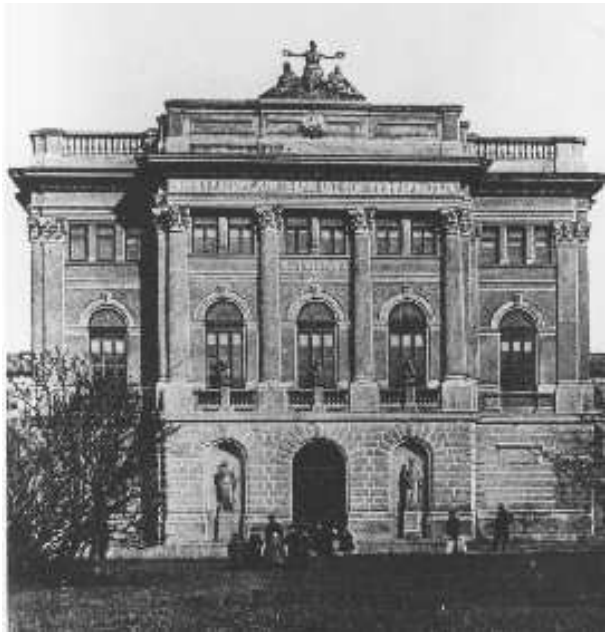}}
\end{center}
\begin{center}
\vskip 0.0truecm
\caption{The University of Warsaw around 1900, at the time Georgy Vorono\"{\i} was professor 
at the university. The photo shows central library building.}
\end{center}
\vskip -1.00truecm
\label{fig:zhuravkawarsaw}
\end{figure*}

Vorono\"{\i} treated his lectures very seriously. He tried to
acquaint his students with recent achievements of science and with
his own new results. To ensure better understanding of the course by
the students, Vorono\"{\i} repeatedly  asked for a permission to
give additional lectures on analytic geometry\footnote{ Though
such supplementary lectures were unpaid for, a lecturer had to
obtain the permission from not only Dean and Rector, but also from
Trustee of the Educational District.}. In the autumn of 1898, Vorono\"{\i} 
also became a Professor at the Warsaw Polytechnical Institute, where 
he became a Dean of the Faculty of Mechanics. In 1898, Vorono\"{\i} got 
elected a member of the Moscow Mathematical Society.

\begin{figure*}[h]
\begin{center}
\mbox{\hskip 0.0truecm\includegraphics[width=11.8cm]{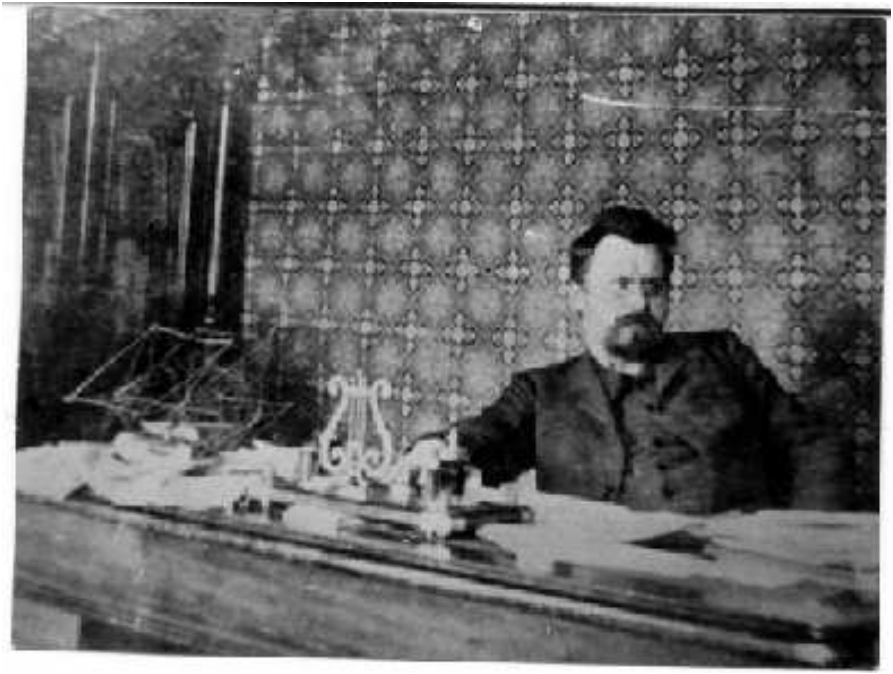}}
\end{center}
\begin{center}
\vskip 0.0truecm
\caption{Georgy Vorono\"{\i} behind his desk. The photo was made around 1902-1903.}
\end{center}
\vskip -1.0truecm
\label{fig:warsawdesk}
\end{figure*}

\noindent In August of 1898, Vorono\"{\i} took part in the tenth Conference of
Russian naturalists and physicians in Kyiv. In 1901, Vorono\"{\i}
was also a participant of the next, eleventh, Conference of Russian
naturalists and physicians in St.~Petersburg and presented three
reports. In one of them, he proposed an original method for
generalized summation of divergent series. Later, in 1919, the
same technique was independently proposed by the Dutch mathematician
N.E. N\"{o}rlund, so that for a long time it was known as N\"{o}rlund's
method. In 1904, Vorono\"{\i} presented two reports at the International
Mathematical Congress held in Heidelberg. Here he met H.~Minkowski, and 
they discovered they were both working on similar topics. 

In 1903 and 1904, Vorono\"{\i} submitted two large papers on Analytic and Algebraic 
Number Theory, {\it ``Sur un probl{\`e}me du calcul des fonctions asymptotiques''} \cite{vor1903} 
and {\it ``Sur une fonction transcendente et ses applications \`a la sommation de quelques 
s\'eries''} \cite{vor1904}. This was a new direction of his studies. In these 
he established the existence of an explicit formula for sums of the form
\begin{equation}
\sum_{n \geq 1}\,d(n)\,f(n)\,\nonumber
\end{equation}
for a large class of functions $f(n)$, including characteristic functions of bounded intervals, and 
$d(n)$ is the number of positive divisors of $n$ \cite{millschmid2003,millschmid2006}. By allowing 
these sums over $f(n)$ to be weighted, the {\it Vorono\"{\i} summations} introduce more general integral 
operations on $f$ than the Fourier transform. These results obtained by Vorono\"{\i} were highly 
appreciated by mathematicians, and as a result of this he was elected Corresponding Member of 
the Russian Academy of Sciences in 1907.

In between 1905 and 1907, both Warsaw University and the Warsaw
Polytechnical Institute were closed  because of revolutionary
events. A group of their staff, amongst whom Vorono\"{\i}, were
sent to Novocherkassk (Southern part of Russia), where the Donskoi
Polytechnical Institute was founded. Vorono\"{\i} worked there
for about a year as a Dean of the Faculty of Mechanics. In autumn 1908,
studies at the Warsaw University were resumed and Vorono\"{\i}
returned to Warsaw. His teaching load was enormous because for 
some time he remained the only professor. Later, Vorono\"{\i}
handed part of his courses to Prof. I.R.~Braitsev, who was transferred 
from Warsaw Polytechnical Institute to assist him. At the time  
Vorono\"{\i} gave a new course on mathematical analysis and wrote 
a textbook on the basis of his lectures. This textbook got published 
in 1909-1911 by the Warsaw University press, of which Braitsev was editor. 
Separately, the book got published in Russian in Kyiv in 1914.  

\begin{figure*}[b]
\begin{center}
\mbox{\hskip 0.0truecm\includegraphics[width=12.0cm]{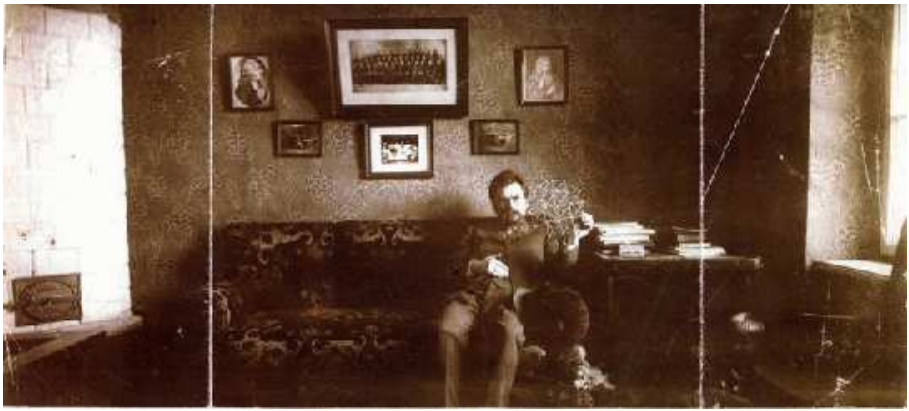}}
\end{center}
\begin{center}
\vskip 0.0truecm
\caption{Georgy Vorono\"{\i} contemplating crystal forms.}
\end{center}
\vskip -1.00truecm
\label{fig:sofa}
\end{figure*}

\section{Quadratic Forms \& Voronoi Tessellations: 1908}
\noindent For many years, Georgy Vorono\"{\i} was working on the development of 
the arithmetic theory of quadratic forms and the geometry of numbers. He had 
the particular ability to think his ideas over and over again until they matured 
and acquired the appropriate perfect form. Once he reached this point, he was able to 
write down the obtained results very fast. 

Preceded by the related study {\it ``Sur quelques propri{\'e}t{\'e}s des formes quadratiques 
positives parfaites''} published in 1908 \cite{vor1908a}, Vorono\"{\i} started to write 
his manuscript on the theory of parallelohedra on March 25, 1907. In ten days he completed 
the work of 106 pages, large size pages filled with a small 
letter font. Yet, barely two weeks later, on April 5, he rewrote 
the text, inserting numerous corrections. Only after the third revision,
Vorono\"{\i} sent his manuscript {\it ``Recherches sur les parall\'elo{\` e}dres primitifs''} 
\cite{vor1908b,vor1909} to the journal, edited by A.L. Krelle, with the following covering letter 
\footnote{The letter (written in French) is kept at the Vorono\"{\i}'s archive of the
Institute of Manuscripts of the National Library of Ukraine.}:

\begin{figure*}[t]
\begin{center}
\vskip 0.0truecm
\framebox[11.8truecm]{\hskip 0.0truecm\includegraphics[width=11.6cm]{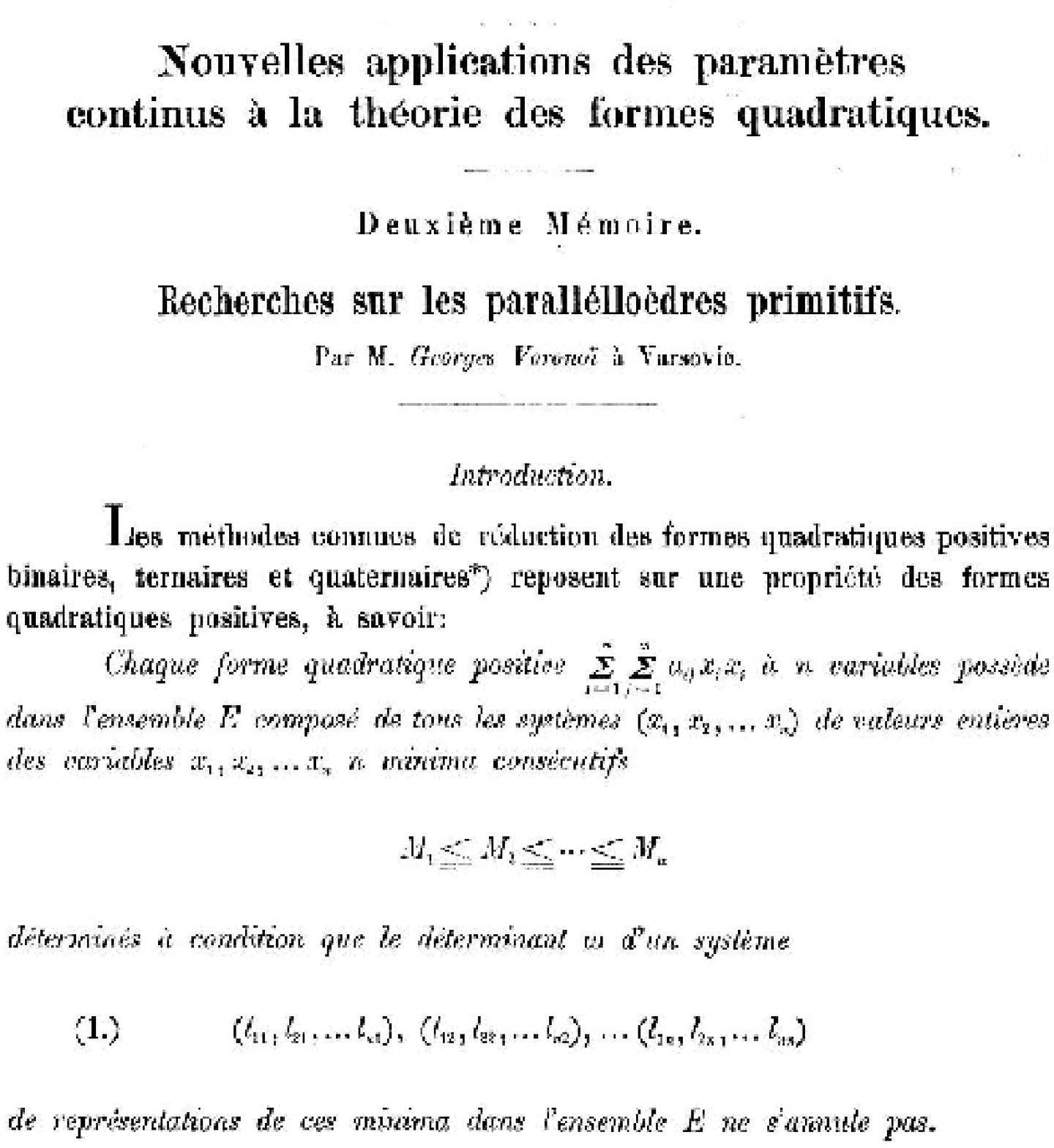}}
\end{center}
\begin{center}
\vskip 0.0truecm
\caption{First page of Vorono\"i's 2nd seminal 1908 paper ... Reproduced with permission 
of Journal f\"ur die reine und angewandte Mathematik.}
\end{center}
\vskip -1.0truecm
\label{fig:paper1908}
\end{figure*}

"For twelve years I have been studying properties of parallelohedra. I can say it 
is a thorny field for investigation, and the results which I obtained and set forth in 
this memoir cost me dear...

Three-dimensional parallelohedra are now playing an important role
in the theory of crystalline bodies, and crystallographers have
already paid attention to properties of these strange polyhedra,
but till now crystallographers were satisfied with the
description of parallelohedra from a purely geometrical point of
view. I noticed already long ago that the task of dividing the
n-dimensional analytical space into convex congruent polyhedra is
closely related to the arithmetic theory of positive quadratic
forms."

Published in 1908-09, this paper, which was certainly the highest manifestation of his
great intelligence, became his Swan Song ... and perhaps the work for which he has 
received greatest acclaim among many different branches of science, Voronoi Tessellations 
(for an example, see fig.~\ref{fig:vor3d}). 

Today Voronoi tessellations have wide applications to the analysis of spatially 
distributed data, ranging from fields such as geophysics, astrophysics and 
cosmology to biology and archaeology. Under different names, such as Dirichlet cells 
and Wigner-Seitz cells the notion is also found in condensed matter physics and 
in the study of Lie groups. 

\section{Failing Health}
\noindent The intense intellectual work took a lot of his energy. 
However, Vorono\"{\i} lacked a good and strong health. Some years 
before physicians had found a disease of the gallbladder which caused 
him a lot of physical and psychological agony. At the time, he worked in 
the new field of indefinite quadratic forms, he spoke about his 
ailments in a rather emotional way. 

>From his mathematical diary we find that Vorono\"{\i} first recorded his 
results on the theory of indefinite quadratic forms on February 20, 1908, 
while in Novocherkassk. Adverse conditions of his life in Novocherkassk 
had significantly worsened his disease. In view of this situation, 
Vorono\"{\i} hurried to write down his thoughts in his mathematical diary.

\smallskip
"I am advancing very well with respect to the problem on which 
I am working; meanwhile my health is getting ever worse. Yesterday for 
the first time I got a clear idea of the algorithm which is to solve all
problems concerning the theory of forms which I am working on, but
just yesterday I had a strong attack of a bilious colic which
prevented me from working in the evening and was keeping me awake
almost all night. I am so afraid lest the results of my long
efforts, which I got with a great deal of hard work, should perish
with me, while it is so difficult to put them in order. Many
things I can only guess by some feeling, which just now, while I am
ill, became sharper ..."
\smallskip

Physicians believed that Georgy Vorono\"{\i} needed a long-term
holiday. But he could not live without his favourite work. Once
Vorono\"{\i} confessed to one of his friends:

\smallskip
"The doctors have forbidden me to work... Moreover, I noticed
myself that a strong intellectual strain always has an effect on
my disease. But they do not know what not being occupied with
mathematics means to me... For me mathematics is life,
everything ..." 
\smallskip

The doctors advised him to go to Carlsbad\footnote{Now it is
Karlovy Vary in the Czech Republic.} for treatment. But he decided to
spend all summer, as in previous years, in Zhuravka, where he had
always restored his health. Indeed, Vorono\"{\i} felt fine but at
the end of October in Warsaw his disease passed on to the acute
form, and Georgy Vorono\"{\i} died on the 7(20)th of November.

\medskip
The early death of Georgy Vorono\"{\i} struck everybody who knew him.
In his obituary professor I.~Braitsev wrote \cite{Br}:
"Nobody could believe that Georgy Theodosi-yovych, who was
highly respected and loved by everybody, had died. Something
extraordinary had happened -- such was the general feeling.
Everybody realized that they had prematurely lost an outstanding
scientist, renowned professor who had been a pride and an
adornment of two higher schools in Warsaw... Seeing off the
deceased's remains to the railway station for carrying them to the
place of burial in the village of Zhuravka, everybody was grieved
the more so because they had forever lost a truthful, sympathetic
and warm-hearted man..."

Georgy Vorono\"{\i}'s body, according to his last will, was
transported to Zhuravka, to his native land. Vorono\"{\i} was
embalmed and interred in a specially built crypt. Later on, in 1910, 
his father, Theodosii Voronyi, was also buried in this crypt. 
In 1932, at the period of the collectivization and the communist 
terror and struggle against the kulaky (the well-to-do peasants) in
the Soviet Union, this crypt was destroyed. The Voronyi family had to leave 
the village that very night -- it was dangerous to remain there. They did 
not return to their native place anymore... The peasants reburied the remains 
of Georgy Vorono\"{\i} and his father in the common grave nearby.

\section{Aftermath: the Vorono\"{\i} family}
\noindent The further destiny of the family was rather sad and grievous 
in the years of the communist great terror \footnote{These data were
communicated to H. Syta by Vorono\"{\i}'s daughter Maria Vorona-Vasylenko 
(1900-1984) and confirmed by other relatives.}.

His wife, Olha, was a midwife. When in Zhuravka she treated peasants free 
of charge\footnote{Her brother, Borys Mytrophanovych Krytskii, was a 
well-known physician in the Pryluky region.}. When the family had to leave 
Zhuravka in 1932, Olha Mytrophanivna went to her eldest son, Oleksander, 
who lived in Yahotyn, a small town near Pryluky. There she lived and died 
in 1939 and was buried on the bank of the local lake the Supoi. Oleksander was 
a physician. He believed that an electrical current in a human body, 
operated in a specific way, could destroy cancer cells. He succeeded in curing  
some kinds of tumours, mainly in the stomach. In 1938 he was arrested 
and did not return from one of Stalin's hard labour camps.

The eldest daughter, Oleksandra, had been educated at the Women's courses of higher 
level, and taught at the special institute for advanced studies of Ukrainian leaders.
\footnote{So-called "faculty of a special purpose".} She also published 
articles on Ukrainian literature. Her husband, Kostyantyn Polubotko, was a
research worker at the Ukrainian Academy of Sciences. In 1937 he
was arrested and exiled to the North, to Ohotske Sea, and had to
work as a fisherman. He had a poor health and within three years 
he died there. Oleksandra died in a house for aged people. 

The youngest son, Yuri, was a surgeon. He was the first person in
the world who accomplished transplantation of a human organ, in 1933 
in Kharkiv. The transplanted kidney was functioning for two days. He could 
not continue his work in the oppressive political atmosphere which prevailed 
in the country. Yuri's wife, Vira Nechayivska, was actively participating in 
Ukrainian political life at the times of the Ukrainian Central Rada headed by 
M.~Hrushevsky. In 1917, she was a member of the government as a representative 
of the Women's Council \cite{Verst}.

The younger daughter, Maria, was a school teacher. Her husband,
Vasyl Vasylenko, was a veterinary surgeon. In 1937, he was
asked by the secret police to become an informer and to inform
about his colleagues. He rejected the proposal and was sentenced
to death. Later on Maria was arrested too and sentenced to ten
years of the camp of special regime in the Far East. The youngest 
daughter, Olena, became a dentist. 

Old people in Zhuravka recall some episodes of the sad 
events during the Holodomor, the Soviet famine of 1932-1933, and 
the Great Terror.  The villagers of Zhuravka remember that the trees 
of a luxuriant garden, which had been grown by Theodosii Voronyi, were
cut. "There were apple trees of many kinds, no more such trees 
grew in Zhuravka thereafter'', one old resident relates, 
``When blossoming they were so fragrant that one could have a headache. 
Voronyis were very diligent people. G.~Vorono\"{\i}'s wife was a small woman 
and very quick. One could always see her hurry along a village street -- to 
provide medical care to people. On Christmas eve they put up a large
Christmas tree for all the children in Zhuravka..."

The Voronyi's house in Zhuravka existed up to 1993, and was used as
a building for the small pupils' school. In 1993, when the new
school building was erected, it began to go to ruin and gradually
was razed to the ground by Zhuravka inhabitants. Now  the street
on which this house was is called "Vorono\"{\i} Street" and the
local school named after G. Vorono\"{\i}. There is a small 
exhibition about the Voronyis family in the local school museum.

\begin{figure*}[h]
\begin{center}
\vskip 0.0truecm
\framebox[11.8truecm]{\hskip 0.0truecm\includegraphics[width=11.5cm]{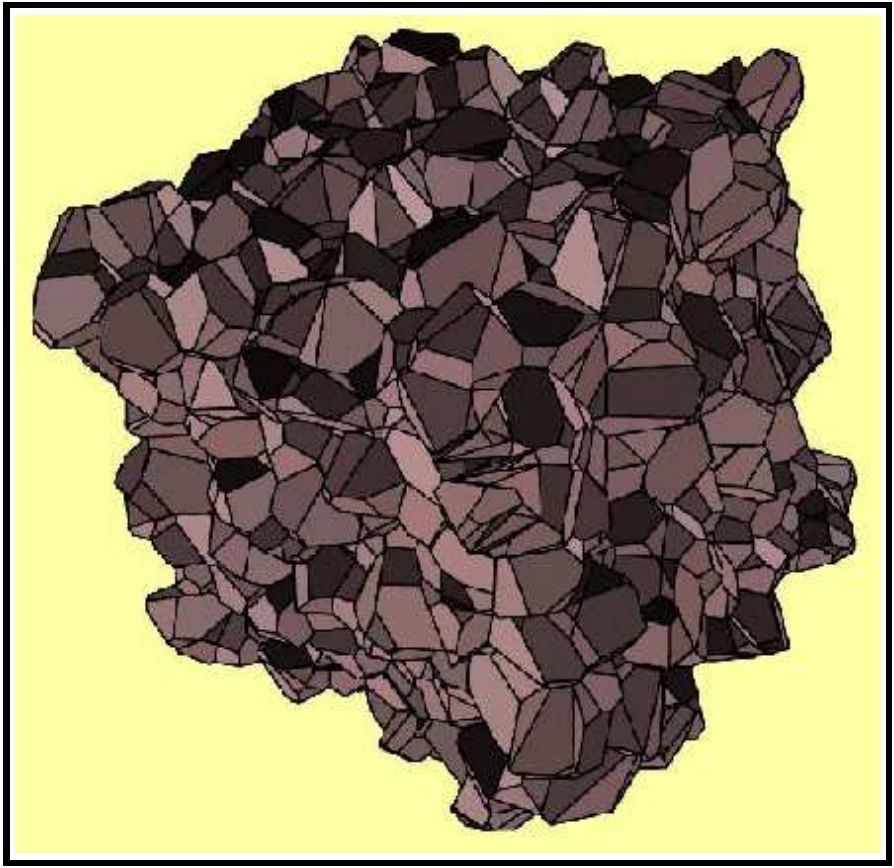}}
\end{center}
\begin{center}
\vskip -0.2truecm
\caption{An example of a 3-D Voronoi tessellation, generated by a set of 1000 Poissonian 
distributed points within a periodic box. Image courtesy of Jacco Dankers.}.
\end{center}
\vskip -1.0truecm
\label{fig:vor3d}
\end{figure*}

\section{Legacy}
Georgy Vorono\"{\i}'s children saved the manuscripts, note-books
of the mathematical diary, and other documents that had belonged to 
him. They transfered them to the possession of the
National Library of Ukraine, to its Institute of Manuscripts where
at present they are kept. There is a description of Vorono\"{\i} archives
in \cite{Pogr,Ven}. In 1952--1953, the Complete Works of 
G.Vorono\"{\i} (in three volumes) were published by the Institute 
of Mathematics of Ukraine with the detailed comments by B.N.~Delone, 
B.A.~Venkov, Yu.V.~Linnyk, I.B.~Pogrebyssky, I.Z.~Shtokalo. 

In the middle of the previous century Boris~Delone, one 
of the well-known followers of Vorono\"{\i}'s research, wrote \cite{Del2}:

"Vorono\"{\i} investigated mainly the problems of Number
Theory. During his short life he published not so many works --
six large memoirs and six short papers. But the profundity and
the significance of his vast investigations left a deep trace in
modern Number Theory. Along with Minkowski, Vorono\"{\i} is
a founder of the Geometry of Numbers.

Vorono\"{\i}'s work of 1903 on the number of points under a hyperbola has to be 
considered as the landmark from which modern Analytic Number Theory begins. 

In Vorono\"{\i}'s work on the algorithm for the calculation of the cubic units, he 
put a series of problems on the distribution of relative minimum. One of these problems, 
amongst the most difficult ones, Vorono\"{\i} has solved himself. Most of them are still
waiting for their solution.''

\bigskip
Following his death, I.~Braitsev expressed his deep regret that
Vorono\"{\i} only succeeded to write down a very small part
of his extensive work on the theory of indefinite quadratic forms 
on paper.

"Hardly can we hope to recreate from them at least in part
those sophisticated trains of geometrical thought which were 
mentioned in the diary algorithm, and of which the deceased always
spoke with great inspiration and enthusiasm. It is not enough to know 
in general outline those directions which the deceased followed. 
One should be such a brilliant expert in the theory of quadratic forms 
with n-variables as he was. It would be 
necessary to wield the marvelous technnique which the deceased
mastered at the end of his life, and it is also necessary to be so
staunchly devoted to this branch of mathematics as he was."

On this issue of the geometry of quadratic forms, Boris Delone 
wrote that ``...Vorono\"{\i}'s investigation on quadratic forms and space 
packing (to which he had added such fundamental results) has not yet exhausted 
the important questions raised by him ... '' \cite{Del2}. Indeed, it made him 
make the following remark to characterize the impact of Vorono\"{\i} short career:

\bigskip
\hskip 0.7truecm\parbox[c]{9.5truecm}{\it
`` Vorono\"{\i}'s memoir on parallelohedra represents one of the deepest 
investigations in the geometry of numbers in the world's literature, and 
the originality of the methods employed in the purely geometrical first part 
stamps the memoir with the imprint of genius.''}

\bigskip
\noindent It shows that the research of G.Vorono\"{\i}'s was already recognized as remarkable by his
contemporaries. In our time, it remains the source of inspiration 
for a large diversity of applications ... 

\bigskip
\noindent {\bf Acknowledgements} The authors are grateful to Vincent Icke and Gert Vegter for 
many helpful suggestions for improving the manuscript. We also wish to acknowledge 
prof.~O.~Ganiyushkin from Kyiv University for advice on algebraic number theory, and 
to Jacco Dankers for producing fig.~\ref{fig:vor3d}.

\end{document}